\documentclass[12pt,reqno]{amsart}

\usepackage{amssymb}

\usepackage{eucal}

\input{diagrams}

\font\sss=cmss8

\def\cG{{\mathcal G}}

\def\BZ{{\mathbb Z}}

\def\sE{\mbox{\sf E}}

\def\sK{\mbox{\sf K}}

\def\ast{{\textstyle *}}

\def\Coker{\operatorname{Coker}}

\def\Ext{\operatorname{Ext}}

\def\GExt{\operatorname{Ext}_{\cG}}

\def\H{\operatorname{H}}

\def\Hom{\operatorname{Hom}}
\def\id{\operatorname{id}}
\def\Image{\operatorname{Im}}

\def\Ker{\operatorname{Ker}}

\def\Mod{\mbox{\sf Mod}}

\def\Proj{\mbox{\sf Pro}}

\def\TExt{\widehat{\operatorname{Ext}}}

\numberwithin{equation}{part}

\newtheorem{Lemma}{Lemma}[section]
\newtheorem{Theorem}[Lemma]{Theorem}
\newtheorem{Proposition}[Lemma]{Proposition}

\theoremstyle{definition}
\newtheorem{Definition}[Lemma]{Definition}
\newtheorem{Setup}[Lemma]{Setup}

\newtheorem{Construction}[Lemma]{Construction}
\newtheorem{Remark}[Lemma]{Remark}

\def\R{A}
\def\matlisR{B}

\def\KProjR{\sK(\Proj\,\R)}
\def\KProjRsmall{\mbox{\sss K}(\mbox{\sss Pro}\,\R)}
\def\EnochsR{\sE(\R)}
\def\EnochsRsmall{\mbox{\sss E}(\R)}

\def\inc{e_{\ast}}
\def\adj{e^!}
\def\res{\operatorname{res}}

\def\Rlm{$\R$-left-module}

\def\M{M}          
\def\N{N}       

\def\Gorprojmod{G}      
\def\kermod{Z}          
\def\projmod{Q}         

\def\kercomplx{K}       
\def\kercomplxprime{K^{\prime}}
\def\projcomplx{P}      

\def\Enochscomplx{E}    

\def\Enochscomplxtil{F} 
\def\Enochscomplxtilprime{F^{\prime}}
\def\Tatecomplx{T}      
\def\Gorprojcomplx{G}   

\def\Enochscomplxmap{e} 
 
\def\Gorprojmodmap{g}   
\def\Tatecomplxmap{t}   

\begin{document}

\title[Tate cohomology]
{Tate cohomology over fairly general rings}

\author{Peter J\o rgensen}
\address{Department of Pure Mathematics, University of Leeds,
Leeds LS2 9JT, United Kingdom}
\email{popjoerg@maths.leeds.ac.uk, www.maths.leeds.ac.uk/\~{ }popjoerg}




\begin{abstract} 

Tate cohomology was originally defined over finite groups.  More
recently, Avramov and Martsinkovsky showed how to extend the
definition so that it now works well over Gorenstein rings.

This paper improves the theory further by giving a new definition that
works over more general rings, specifically, those with a dualizing
complex.

The new definition of Tate cohomology retains the desirable properties
established by Avramov and Martsinkovsky.  Notably, there is a long
exact sequence connecting it to ordinary Ext groups.

\end{abstract}

\maketitle

\setcounter{section}{-1}
\section{Introduction}
\label{sec:introduction}

Tate cohomology was originally defined over finite gro\-ups, and has
been used to great effect in group representation theory.

More recently, Avramov and Martsinkovsky accomplished in
\cite{AvrMart} an extension of the definition so that it now works
well over Gorenstein rings.  In fact, \cite{AvrMart} went so far as
to define Tate $\Ext$ groups,
\begin{equation}
\label{equ:Tate_Exts}
  \TExt_{\R}^i(\M,\N),
\end{equation}
which have classical Tate cohomology as the special case
$\TExt_{\BZ G}^i(\BZ,\N)$.

Moreover, a new key result was introduced in \cite{AvrMart}: Tate and
ordinary $\Ext$ groups fit into a long exact sequence
\begin{equation}
\label{equ:long_exact_sequence}
 0 \rightarrow \GExt^1(\M,\N) 
 \longrightarrow \Ext^1(\M,\N)
 \longrightarrow \TExt^1(\M,\N) 
 \longrightarrow \cdots,
\end{equation}
where the $\GExt$'s are relative $\Ext$ groups defined by means of
what is known as (proper) Gorenstein projective resolutions in the
first variable.  This illuminates the connection between Tate
cohomology and ordinary $\Ext$ groups, and also explains why
``Gorenstein phenomena'' play an important role in Tate cohomology
theory.

To be more precise about \cite{AvrMart}, what it really did was to
construct the Tate $\Ext$ groups from equation \eqref{equ:Tate_Exts}
over any noetherian ring, but only when $\M$ is a finitely generated
module of finite Gorenstein projective dimension.  So to have the Tate
$\Ext$ groups defined everywhere, it is necessary to require that each
finitely generated module has finite Gorenstein projective dimension.
In turn, this is the same as requiring the ring to be Gorenstein.
Hence, for non-Gorenstein rings there are modules where the Tate
$\Ext$ groups of \cite{AvrMart} remain undefined.

The present paper improves the theory by giving a new definition of
Tate $\Ext$ groups that works for all modules, including the ones with
infinite Gorenstein projective dimension.  Moreover, it is proved that
the new definition coincides with the definition from \cite{AvrMart}
when they both work, and that the new definition retains the desirable
properties established in \cite{AvrMart}, notably the long exact
sequence \eqref{equ:long_exact_sequence}.

The results apply over fairly general rings, whence the title of the
paper.  For instance, noetherian commutative rings with dualizing
complexes are covered, and so too is a large class of non-commutative
rings.  See remark \ref{rmk:blanket} for more information.

Let me close the introduction with a synopsis of the paper.

Section \ref{sec:setup} gives the setup which is used in the rest
of the paper.

Section \ref{sec:Tate} gives the new definition of Tate $\Ext$ groups
(definition \ref{def:Tate}), proves that short exact sequences of
modules result in long exact sequences of Tate $\Ext$ groups
(proposition \ref{pro:long_exact_sequences}), recalls the definition
of Tate $\Ext$ groups which was given in \cite{AvrMart} (remark
\ref{rmk:AvrMart_Tate}), and proves that the definitions coincide when
they both work (proposition \ref{pro:correspondance}).

Sections \ref{sec:lemmas} and \ref{sec:resolution} form an interlude.
Section \ref{sec:lemmas} just contains some lemmas.  Section
\ref{sec:resolution} considers Gorenstein projective modules, and
shows how the machinery of the paper gives rise to Gorenstein
projective resolutions (lemma \ref{lem:resolution}).  

Section \ref{sec:sequence} proves that the new definition of Tate
$\Ext$ groups fits into the long exact sequence
\eqref{equ:long_exact_sequence} (theorem
\ref{thm:long_exact_sequence}).

As a spin off along the way, it is proved that each module $\M$ admits
a Gorenstein projective precover $\Gorprojmod \longrightarrow \M$
whose kernel $\kermod$ is particularly nice.  Namely, $\kermod$ has
the property that each projective resolution of $\kermod$ is even a
(proper) Gorenstein projective resolution of $\kermod$ (theorem
\ref{thm:special_precovers}).  This generalizes previous results by
several authors (see remark \ref{rmk:special_precovers}).

\setcounter{section}{0}
\section{Setup}
\label{sec:setup}

\begin{Definition}
\label{def:Enochs}
Let $\R$ be a ring.  By $\EnochsR$ is denoted the class of complexes
$\Enochscomplx$ of \Rlm s so that $\Enochscomplx$ consists of
projective modules, is exact, and has
$\Hom_{\R}(\Enochscomplx,\projmod)$ exact for each projective \Rlm\
$\projmod$.
\end{Definition}

I will view $\EnochsR$ as a full subcategory of $\KProjR$, the
homotopy category of complexes of projective \Rlm s.

Note that $\EnochsR$ consists precisely of the complexes known as
complete projective resolutions.

\begin{Setup}
\label{set:blanket}
Throughout this paper, $\R$ is a ring for which the inclusion functor
\[
  \inc : \EnochsR \longrightarrow \KProjR
\]
has a right-adjoint
\[
  \adj : \KProjR \longrightarrow \EnochsR.
\]
\end{Setup}

\begin{Remark}
\label{rmk:blanket}
The right-adjoint $\adj$ exists for fairly general rings, whence the
title of the paper.

Namely, by \cite[prop.\ 2.2]{PJgorproj} the right-adjoint $\adj$
exists when $\R$ is a noetherian commutative ring with a dualizing
complex.

Also, by \cite[sec.\ 4]{PJgorproj} the right-adjoint $\adj$
exists when $\R$ is a left-coherent and right-noetherian $k$-algebra over
the field $k$ for which there exists a left-noetherian $k$-algebra
$\matlisR$ and a dualizing complex ${}_{\matlisR}D_{\R}$.

These cases cover many rings arising in practice.
\end{Remark}

\begin{Remark}
\label{rmk:approximation}
If $\projcomplx$ is a complex of projective modules, then $\adj
\projcomplx$ can be thought of as the best approximation to
$\projcomplx$ by a complete projective resolution.

Elaborating on this, if $\M$ is a module with projective
resolution $\projcomplx$, then $\adj \projcomplx$ can be thought of as
the best approximation to $\M$ by a complete projective
resolution.  This point will be made more precise in lemma
\ref{lem:adj_gives_complete_resolutions}.
\end{Remark}

\section{Tate $\Ext$ groups}
\label{sec:Tate}

This section gives the new definition of Tate $\Ext$ groups
(definition \ref{def:Tate}), and proves that short exact sequences of
modules result in long exact sequences of Tate $\Ext$ groups
(proposition \ref{pro:long_exact_sequences}).

The rest of the section is devoted to recalling the earlier definition
of Tate $\Ext$ groups which was given in \cite{AvrMart} (remark
\ref{rmk:AvrMart_Tate}), and proving that the new and the earlier
definition of Tate $\Ext$ groups coincide when
they both work (proposition \ref{pro:correspondance}).

\begin{Remark}
It is classical that the category of \Rlm s $\Mod(\R)$ is equivalent
to the full subcategory of $\KProjR$ consisting of projective
resolutions of \Rlm s.  Let 
\[
  \Mod(\R) \stackrel{\res}{\longrightarrow} \KProjR
\]
be a functor implementing the equivalence.
\end{Remark}

\begin{Definition}
\label{def:Tate}
For \Rlm s $\M$ and $\N$, the $i$'th Tate Ext group is
\[
  \TExt^i(\M,\N) = \H^i\!\Hom_{\R}(\adj \res \M,\N).
\]
\end{Definition}

This is the definition one would expect: As pointed out in remark
\ref{rmk:approximation}, the complex $\adj \res \M$ can be
thought of as the best approximation to $\M$ by a complete
projective resolution.  So taking $\Hom$ into $\N$ and taking
cohomology should be the way to get Tate $\Ext$ groups.

\begin{Proposition}
\label{pro:long_exact_sequences}
If 
\[
  0 \rightarrow \M^{\prime} 
    \longrightarrow \M 
    \longrightarrow \M^{\prime \prime} \rightarrow 0 
\]
and 
\[
  0 \rightarrow \N^{\prime}
    \longrightarrow \N 
    \longrightarrow \N^{\prime \prime} \rightarrow 0
\]
are short exact sequences of \Rlm s, then there are natural long exact
sequences
\[
  \cdots \longrightarrow \TExt^i(\M^{\prime \prime},\N)
  \longrightarrow \TExt^i(\M,\N)
  \longrightarrow \TExt^i(\M^{\prime},\N)
  \longrightarrow \cdots
\]
and
\[
  \cdots \longrightarrow \TExt^i(\M,\N^{\prime})
  \longrightarrow \TExt^i(\M,\N)
  \longrightarrow \TExt^i(\M,\N^{\prime \prime})
  \longrightarrow \cdots.
\]
\end{Proposition}

\begin{proof}
It is well known that the first short exact sequence in the
proposition results in a distinguished triangle in $\KProjR$,
\[
  \res \M^{\prime} 
  \longrightarrow \res \M
  \longrightarrow \res \M^{\prime \prime}
  \longrightarrow.
\]
Since $\inc$ is a triangulated functor, so is its adjoint $\adj$, so
there is also a distinguished triangle in $\EnochsR$,
\[
  \adj \res \M^{\prime} 
  \longrightarrow \adj \res \M
  \longrightarrow \adj \res \M^{\prime \prime}
  \longrightarrow.
\]
This again results in a distinguished triangle
\[
  {\scriptstyle
  \Hom_{\R}(\adj \res \M^{\prime \prime},\N)
  \longrightarrow \Hom_{\R}(\adj \res \M,\N)
  \longrightarrow \Hom_{\R}(\adj \res \M^{\prime},\N)
  \longrightarrow
  }
\]
whose cohomology long exact sequence is the first long exact sequence
in the proposition. 

The complex $\adj \res \M$ is in $\EnochsR$ so consists of projective
modules, so the second short exact sequence in the proposition gives a
short exact sequence of complexes
\[
  {\scriptstyle
  0 \rightarrow \Hom_{\R}(\adj \res \M,\N^{\prime})
    \longrightarrow \Hom_{\R}(\adj \res \M,\N)
    \longrightarrow \Hom_{\R}(\adj \res \M,\N^{\prime \prime})
    \rightarrow 0
  }
\]
whose cohomology long exact sequence is the second long exact sequence
in the proposition. 
\end{proof}

\begin{Construction}
\label{con:surjective}
If $\projcomplx$ is a complex of \Rlm s, then for each $i$ there is a
chain map
\[
  \begin{diagram}[labelstyle=\scriptstyle,width=4ex,height=7ex]
    \cdots & \rTo & 0 & & \rTo & & \projcomplx^i & & \rTo^{\id} & & \projcomplx^i & & \rTo & & 0 & \rTo & \cdots \\
    & & \dTo & & & & \dTo^{\id} & & & & \dTo_{\partial_{\projcomplx}^i} & & & & \dTo & & \\
    \cdots & \rTo & \projcomplx^{i-1} & & \rTo_{\partial_{\projcomplx}^{i-1}} & & \projcomplx^i & & \rTo_{\partial_{\projcomplx}^i} & & \projcomplx^{i+1} & & \rTo_{\partial_{\projcomplx}^{i+1}} & & \projcomplx^{i+2} & \rTo & \cdots \\
  \end{diagram}
\]
where the upper complex is null homotopic.  

This is useful because I can add the upper complex to $\Tatecomplx$ in
any chain map $\Tatecomplx \stackrel{\Tatecomplxmap}{\longrightarrow}
\projcomplx$, and thereby change $\Tatecomplxmap$ so that the $i$'th
component $\Tatecomplxmap^i$ becomes surjective.  Doing so does not
change the isomorphism class of $\Tatecomplxmap$ in $\sK(\R)$, the
homotopy category of complexes of \Rlm s.
\end{Construction}

\begin{Remark}
\label{rmk:AvrMart_Tate}
The earlier definition of Tate $\Ext$ groups which was given in
\cite{AvrMart} is
\[
  \TExt^i(\M,\N) = \H^i\!\Hom_{\R}(\Tatecomplx,\N)
\]
where $\Tatecomplx$ is a complete projective resolution of the \Rlm\
$\M$.  This means that $\Tatecomplx$ is in $\EnochsR$,
consists of finitely generated modules, and sits in a diagram of chain
maps
\begin{equation}
\label{equ:complete_resolution}
  \Tatecomplx \stackrel{\Tatecomplxmap}{\longrightarrow} 
  \projcomplx \longrightarrow \M
\end{equation}
where $\projcomplx \longrightarrow \M$ is a projective resolution and
where $\Tatecomplxmap^i$ is bijective for $i \ll 0$.

Note that not all \Rlm s have complete projective resolutions.  In
fact, if $\R$ is left-noetherian, then the ones that do are exactly
the ones which are finitely generated and have finite Gorenstein
projective dimension by \cite[thm.\ 3.1]{AvrMart}. 
\end{Remark}

\begin{Lemma}
\label{lem:adj_gives_complete_resolutions}
Let $\M$ be an \Rlm\ which has a projective resolution $\projcomplx$
and a complete projective resolution $\Tatecomplx$.  Then
\[
  \adj \projcomplx \cong \Tatecomplx  
\]
in $\KProjR$.
\end{Lemma}

\begin{proof} 
All projective resolutions of $\M$ are isomorphic in $\KProjR$,
so I may as well prove the lemma for the specific projective
resolution $\projcomplx$ from equation
\eqref{equ:complete_resolution}.

By applying construction \ref{con:surjective} to the chain map
$\Tatecomplx \stackrel{\Tatecomplxmap}{\longrightarrow}
\projcomplx$ in cohomological degrees larger than some number, I can
assume that $\Tatecomplxmap$ is surjective.  Hence there is a short
exact sequence of complexes
\begin{equation}
\label{equ:pre_basic_exact_sequence}
  0 \rightarrow \kercomplx 
    \longrightarrow \Tatecomplx
    \stackrel{\Tatecomplxmap}{\longrightarrow} \projcomplx
    \rightarrow 0.
\end{equation}

Since both $\Tatecomplx$ and $\projcomplx$ consist of projective
modules, the sequence is semi-split and $\kercomplx$ also consists of
projective modules.  Moreover, by assumption, $\Tatecomplxmap^i$ is
bijective for $i \ll 0$, so $K^i = 0$ for $i \ll 0$.  So
$\kercomplx$ is a left-bounded complex of projective modules.

Now let $\Enochscomplx$ be in $\EnochsR$.  In particular,
$\Hom_{\R}(\Enochscomplx,\projmod)$ is exact when $\projmod$ is a
projective module.  It is classical that
$\Hom_{\R}(\Enochscomplx,\kercomplx)$ is then also exact, because
$\kercomplx$ is a left-bounded complex of projective modules.  Indeed,
this follows by an argument analogous to the one which shows that if
$X$ is an exact complex and $I$ is a left-bounded complex of injective
modules, then $\Hom_{\R}(X,I)$ is exact.

Since the sequence \eqref{equ:pre_basic_exact_sequence} is semi-split,
it stays exact under the functor $\Hom_{\R}(\Enochscomplx,-)$.  So
there is a short exact sequence of complexes
\[
  0 \rightarrow \Hom_{\R}(\Enochscomplx,\kercomplx)
    \longrightarrow \Hom_{\R}(\Enochscomplx,\Tatecomplx)
    \longrightarrow \Hom_{\R}(\Enochscomplx,\projcomplx)
    \rightarrow 0.
\]
Since $\Hom_{\R}(\Enochscomplx,\kercomplx)$ is exact, the
cohomology long exact sequence shows that there is an isomorphism
\[
  \H^0\!\Hom_{\R}(\Enochscomplx,\Tatecomplx)
  \cong \H^0\!\Hom_{\R}(\Enochscomplx,\projcomplx)
\]
which is natural in $\Enochscomplx$.  That is, there is a natural
isomorphism
\[
  \Hom_{\KProjRsmall}(\Enochscomplx,\Tatecomplx)
  \cong \Hom_{\KProjRsmall}(\Enochscomplx,\projcomplx)
\]
which can also be written 
\[
  \Hom_{\EnochsRsmall}(\Enochscomplx,\Tatecomplx)
  \cong \Hom_{\KProjRsmall}(\Enochscomplx,\projcomplx)
\]
because $\Enochscomplx$ and $\Tatecomplx$ are in $\EnochsR$.

On the other hand, I also have a natural isomorphism
\[
  \Hom_{\KProjRsmall}(\Enochscomplx,\projcomplx)
  = \Hom_{\KProjRsmall}(\inc \Enochscomplx,\projcomplx)
  \cong \Hom_{\EnochsRsmall}(\Enochscomplx,\adj \projcomplx).
\]
Combining the last two equations gives a natural isomorphism
\[
  \Hom_{\EnochsRsmall}(\Enochscomplx,\Tatecomplx)
  \cong \Hom_{\EnochsRsmall}(\Enochscomplx,\adj \projcomplx),
\]
proving $\Tatecomplx \cong \adj \projcomplx$ as desired.
\end{proof}

\begin{Proposition}
\label{pro:correspondance}
Let $\M$ be an \Rlm\ which has a complete projective resolution
$\Tatecomplx$.  Then the new Tate $\Ext$ groups of this paper
(definition \ref{def:Tate}) coincide with the Tate $\Ext$ groups which
were defined in \cite{AvrMart} (remark \ref{rmk:AvrMart_Tate}).
\end{Proposition}

\begin{proof}
Lemma \ref{lem:adj_gives_complete_resolutions} gives that the
projective resolution $\res \M$ of $\M$ satisfies $\adj \res \M \cong
\Tatecomplx$.  Combining this with the formulae in definition
\ref{def:Tate} and remark \ref{rmk:AvrMart_Tate} proves the
proposition.  
\end{proof}

\section{Some lemmas}
\label{sec:lemmas}

This section collects three lemmas needed later in the paper.

\begin{Lemma}
\label{lem:isomorphism}

Let $\projcomplx$ be in $\KProjR$, and consider the counit morphism
\[
  \inc \adj \projcomplx
  \stackrel{\epsilon_{\projcomplx}}{\longrightarrow} 
  \projcomplx.  
\]
Let $\Enochscomplx$ be in $\EnochsR$.  Then the induced map
\[
  \begin{diagram}[labelstyle=\scriptstyle,width=17.5ex,midshaft]
    \Hom_{\KProjRsmall}(\Enochscomplx,\inc \adj \projcomplx)
    & \rTo^{\Hom(\Enochscomplx,\epsilon_{\projcomplx})} &
    \Hom_{\KProjRsmall}(\Enochscomplx,\projcomplx) \\
  \end{diagram}
\]
is an isomorphism.
\end{Lemma}

\begin{proof}
There is a commutative diagram
\[
  \begin{diagram}[labelstyle=\scriptstyle,width=12ex,height=5ex]
    \Hom_{\EnochsRsmall}(\Enochscomplx,\adj \projcomplx)
      & \rTo^{\inc(-)}
      & \Hom_{\KProjRsmall}(\inc \Enochscomplx,\inc \adj \projcomplx) \\
    & \SE & \dTo_{\Hom(\inc \Enochscomplx,\epsilon_{\projcomplx})} \\
    & & \Hom_{\KProjRsmall}(\inc \Enochscomplx,\projcomplx) \lefteqn{.} \\
  \end{diagram}
\]
The diagonal map is the adjunction isomorphism, and the horizontal map
is an isomorphism because $\inc$ is the inclusion functor of a full
subcategory.

Hence the vertical map must be an isomorphism, but this is the induced
map from the lemma in slightly different notation.
\end{proof}

\begin{Remark}
\label{rmk:Gorenstein_projective_epic}
For the following two lemmas, recall that an \Rlm\ is called
Gorenstein projective if it has the form 
\[
  \Gorprojmod = \Ker(\Enochscomplx^1 \longrightarrow \Enochscomplx^2) 
\]
for some $\Enochscomplx$ in $\EnochsR$.

A homomorphism $\kercomplx \stackrel{s}{\longrightarrow} \N$ is called
a relative epimorphism with respect to the class of Gorenstein
projective modules if each homomorphism $\Gorprojmod
\stackrel{\Gorprojmodmap}{\longrightarrow} \N$ with
$\Gorprojmod$ Gorenstein projective lifts through $s$,
\[
  \begin{diagram}[labelstyle=\scriptstyle]
                &                       & \kercomplx \\
                & \ruDotsto             & \dTo_{s} \\
    \Gorprojmod & \rTo_{\Gorprojmodmap} & \N \lefteqn{.}
  \end{diagram}
\smallskip
\]
\end{Remark}

\begin{Lemma}
\label{lem:Gorenstein_projective_epic_1}
Let $\kercomplx$ be a complex of projective \Rlm s satisfying
$\Hom_{\KProjRsmall}(\Enochscomplx,\kercomplx) = 0$ for each
$\Enochscomplx$ in $\EnochsR$, and suppose that
\[
  \cdots \longrightarrow \kercomplx^{i-1} \longrightarrow \kercomplx^i
  \stackrel{s}{\longrightarrow} \N \rightarrow 0
\]
is exact.

Then $\kercomplx^i \stackrel{s}{\longrightarrow} \N$ is a
relative epimorphism with respect to the class of Gorenstein
projective modules.
\end{Lemma}

\begin{proof}
Let $\Gorprojmod$ be a Gorenstein projective module, and let
$\Gorprojmod \stackrel{\Gorprojmodmap}{\longrightarrow} \N$ be
a homomorphism.

By shifting, I can clearly pick a complex $\Enochscomplx$ in
$\EnochsR$ with $G = \Ker(\Enochscomplx^{i+1} \longrightarrow
\Enochscomplx^{i+2})$, and it is not hard to see that there is a chain
map $\Enochscomplx \stackrel{\Enochscomplxmap}{\longrightarrow}
\kercomplx$ which fits together with $\Gorprojmod
\stackrel{\Gorprojmodmap}{\longrightarrow} \N$ in a commutative
diagram 
\[
  \begin{diagram}[labelstyle=\scriptstyle,height=4ex,width=4ex]
    \cdots & \rTo & \Enochscomplx^{i-1} & & \rTo & & \Enochscomplx^i & & \rTo & & \Enochscomplx^{i+1} & & \rTo & & \Enochscomplx^{i+2} & \rTo & \cdots \\
    & & & & & & & \rdOnto & & \ruEmbed_{\ell} & & & & & & & \\
    & & \dTo^{\Enochscomplxmap^{i-1}} & & & & \dTo^{\Enochscomplxmap^i} & & \Gorprojmod & & \dTo_{\Enochscomplxmap^{i+1}} & & & & \dTo_{\Enochscomplxmap^{i+2}} & & \\
    & & & & & & & & \vLine^{\Gorprojmodmap} & & & & & & & & \\
    \cdots & \rTo & \kercomplx^{i-1} & & \rTo & & \kercomplx^i & \rTo & \HonV & & \kercomplx^{i+1} & & \rTo & & \kercomplx^{i+2} & \rTo & \cdots \lefteqn{.} \\
    & & & & & & & \rdOnto_s & \dTo & \ruTo & & & & & & & \\
    & & & & & & & & \N & & & & & & & & \\
  \end{diagram}
\]
Since I have assumed $\Hom_{\KProjRsmall}(\Enochscomplx,\kercomplx) =
0$ for $\Enochscomplx$ in $\EnochsR$, the chain map $\Enochscomplxmap$
must be null homotopic.  

Let $h$ be a null homotopy with $\Enochscomplxmap =
h\partial^{\Enochscomplx} + \partial^{\kercomplx}h$, consisting of
components $\Enochscomplx^j \stackrel{h^j}{\longrightarrow}
\kercomplx^{j-1}$.  Then it is straightforward to prove
\[
  s \circ (h^{i+1} \ell) = \Gorprojmodmap,
\]
so $\Gorprojmod \stackrel{\Gorprojmodmap}{\longrightarrow} \N$
has been lifted through $\kercomplx^i \stackrel{s}{\longrightarrow}
\N$ as desired.
\end{proof}

\begin{Lemma}
\label{lem:Gorenstein_projective_epic_2}
Let $\kercomplx$ be a complex of projective \Rlm s satisfying
$\Hom_{\KProjRsmall}(\Enochscomplx,\kercomplx) = 0$ for each
$\Enochscomplx$ in $\EnochsR$, and suppose that
\[
  \cdots \longrightarrow \kercomplx^{i-2} \longrightarrow \kercomplx^{i-1}
  \longrightarrow \Ker \partial_{\kercomplx}^i
  \stackrel{t}{\longrightarrow} \N \rightarrow 0
\]
is exact.

Then $\Ker \partial_{\kercomplx}^i \stackrel{t}{\longrightarrow}
\N$ is a relative epimorphism with respect to the class of
Gorenstein projective modules.
\end{Lemma}

\begin{proof}
From the data given I can construct a commutative diagram
\[
  \begin{diagram}[labelstyle=\scriptstyle,height=4ex,width=4ex]
    \cdots & \rTo & \kercomplx^{i-1} & & \rTo & & \kercomplx^i & & \rTo & & \kercomplx^{i+1} & \rTo & \cdots \lefteqn{.} \\
    & & & \rdTo & & \ruEmbed & & \rdOnto^s & & \ruTo & & & \\
    & & & & \Ker \partial_{\kercomplx}^i & & & & \Coker \partial_{\kercomplx}^{i-1} & & & & \\
    & & & & & \rdOnto_t & & \ruEmbed & & & & & \\
    & & & & & & \N & & & & & & \\
  \end{diagram}
\]
It follows from lemma \ref{lem:Gorenstein_projective_epic_1} that
$\kercomplx^i \stackrel{s}{\longrightarrow} \Coker
\partial_{\kercomplx}^{i-1}$ is a relative epimorphism with respect to
the class of Gorenstein projective modules, and it is a small diagram
exercise to see that this implies the same for $\Ker
\partial_{\kercomplx}^i \stackrel{t}{\longrightarrow} \N$.
\end{proof}

\section{Gorenstein projective modules}
\label{sec:resolution}

Let $\M$ be an \Rlm\ with projective resolution $\projcomplx$.  This
section considers the kernel $\kercomplx$ of the counit morphism $\inc
\adj \projcomplx \stackrel{\epsilon_{\projcomplx}}{\longrightarrow}
\projcomplx$ (construction \ref{con:basic_exact_sequence}) and shows
how it gives rise to a Gorenstein projective resolution of $\M$ (lemma
\ref{lem:resolution}).

This resolution is constructed to be used in the next section.
However, as a spin off I also use it to prove (theorem
\ref{thm:special_precovers}) that each module $\M$ admits a Gorenstein
projective precover $\Gorprojmod \longrightarrow \M$ whose kernel
$\kermod$ is particularly nice.  Namely, $\kermod$ has the property
that each projective resolution of $\kermod$ is even a Gorenstein
projective resolution of $\kermod$.  This generalizes previous results
by several authors (see remark \ref{rmk:special_precovers}).

\begin{Construction}
\label{con:basic_exact_sequence}

Let $\M$ be an \Rlm\ with projective resolution $\projcomplx$
concentrated in non-positive cohomological degrees.  By applying
construction \ref{con:surjective} in each degree, I can assume that
the counit morphism $\inc \adj \projcomplx
\stackrel{\epsilon_{\projcomplx}}{\longrightarrow}
\projcomplx$ in $\KProjR$ is represented by a surjective chain
map, so setting 
\[
  \Enochscomplxtil = \inc \adj \projcomplx, 
\]
there is a short exact sequence of complexes
\[
  0 \rightarrow \kercomplx 
    \longrightarrow \Enochscomplxtil
    \longrightarrow \projcomplx \rightarrow 0.
\]
Note that since both $\Enochscomplxtil$ and $\projcomplx$ consist of
projective modules, the sequence is semi-split and $\kercomplx$ also
consists of projective modules.
\end{Construction}

\begin{Lemma}
\label{lem:zero}
Consider the situation of construction \ref{con:basic_exact_sequence}.
The complex $\kercomplx$ satisfies
\[
  \Hom_{\KProjRsmall}(\Enochscomplx,\kercomplx) = 0
\]
for $\Enochscomplx$ in $\EnochsR$.
\end{Lemma}

\begin{proof}
The short exact sequence from construction
\ref{con:basic_exact_sequence} is semi-split and therefore gives a
distinguished triangle
\[
  \kercomplx 
  \longrightarrow \Enochscomplxtil
  \longrightarrow \projcomplx 
  \longrightarrow
\]
in $\KProjR$.  Hence there is a cohomology long exact sequence
consisting of pieces 
\[
  {\scriptstyle 
  \Hom_{\KProjRsmall}(\Sigma^i \Enochscomplx,\kercomplx)
  \longrightarrow 
  \Hom_{\KProjRsmall}(\Sigma^i \Enochscomplx,\Enochscomplxtil)
  \longrightarrow
  \Hom_{\KProjRsmall}(\Sigma^i \Enochscomplx,\projcomplx).
  }
\]
Lemma \ref{lem:isomorphism} says that the second map here is always an
isomorphism, and this implies the lemma.
\end{proof}

\begin{Remark}
\label{rmk:Gorenstein_projective_resolutions}
For the lemma below, recall that an augmented Gorenstein projective
resolution of a module $\M$ is an exact sequence
\[
  \cdots \longrightarrow \Gorprojcomplx^{-1} 
  \longrightarrow \Gorprojcomplx^0 
  \longrightarrow \M
  \rightarrow 0
\]
satisfying
\begin{enumerate}

  \item  The modules $\Gorprojcomplx^0, \Gorprojcomplx^{-1}, \ldots$
         are Gorenstein projective.   

  \item  If the sequence is split up into short exact sequences, then
         in each short exact sequence the surjection is a relative
         epimorphism with respect to the class of Gorenstein
         projective modules.

\end{enumerate}

The complex
\[
  \Gorprojcomplx = 
  \cdots \longrightarrow \Gorprojcomplx^{-1} 
  \longrightarrow \Gorprojcomplx^0 
  \longrightarrow 0
  \longrightarrow \cdots
\]
is then called a Gorenstein projective resolution of $\M$.  Note that
some authors call this a {\em proper} Gorenstein projective
resolution.
\end{Remark}

\begin{Remark}
\label{rmk:resolution}
Consider the short exact sequence from construction
\ref{con:basic_exact_sequence}.  The complex $\Enochscomplxtil$ is in
$\EnochsR$.  In particular it is exact, and therefore the cohomology
long exact sequence shows
\[
  \H^i\!\kercomplx =
  \left\{
    \begin{array}{cl}
      \M & \mbox{ for } i =     1, \\
      0       & \mbox{ for } i \not= 1.
    \end{array}
  \right.
\]
Hence there is an exact sequence
\[
  \cdots 
  \longrightarrow \kercomplx^{-2}
  \longrightarrow \kercomplx^{-1}
  \longrightarrow \kercomplx^0
  \longrightarrow \Ker \partial_{\kercomplx}^1
  \stackrel{u}{\longrightarrow} \M
  \rightarrow 0.
\]
\end{Remark}

\begin{Lemma}
\label{lem:resolution}
Let $\M$ be an \Rlm, and consider the exact sequence from remark
\ref{rmk:resolution}.  This is an augmented Gorenstein projective
resolution of $\M$.  
\end{Lemma}

\begin{proof}
I must prove that the sequence satisfies (i) and (ii) from
remark \ref{rmk:Gorenstein_projective_resolutions}.

\medskip
\noindent
(i)  The modules $\kercomplx^0, \kercomplx^{-1}, \ldots$ are
projective and hence Gorensten projective.

As for $\Ker \partial_{\kercomplx}^1$, observe that in the short exact
sequence from construction \ref{con:basic_exact_sequence}, the complex
$\projcomplx$ is concentrated in non-positive cohomological degrees,
so the modules $P^1$ and $P^2$ are zero.  So in degrees $1$ and $2$,
the short exact sequence gives
\[
  \begin{diagram}[labelstyle=\scriptstyle,midshaft]
    \kercomplx^2 & \rTo^{\cong} & \Enochscomplxtil^2 \\
    \uTo^{\partial_{\kercomplx}^1} & & \uTo_{\partial_{\Enochscomplxtil}^1} \\
    \kercomplx^1 & \rTo^{\cong} & \Enochscomplxtil^1 \lefteqn{.} \\
  \end{diagram}
\]
Hence
\[
  \Ker \partial_{\kercomplx}^1 
  \cong \Ker \partial_{\Enochscomplxtil}^1,
\]
and this last module is Gorenstein projective because the complex
$\Enochscomplxtil$ is in $\EnochsR$. 

\medskip
\noindent
(ii)  When splitting up the exact sequence from remark
\ref{rmk:resolution} into short exact sequences, each sequence but
one is equal to a short exact sequence which results from splitting up
$\kercomplx$ itself.  Consider such a sequence
\[
  0 \rightarrow \Ker \partial_{\kercomplx}^i
    \longrightarrow \kercomplx^i 
    \stackrel{k^i}{\longrightarrow} \Image \partial_{\kercomplx}^i 
    \rightarrow 0
\]
with $i \leq 0$.  

Lemma \ref{lem:zero} gives that lemma
\ref{lem:Gorenstein_projective_epic_1} applies to this situation if I
set the homomorphism $\kercomplx^i \stackrel{s}{\longrightarrow} \N$
from lemma \ref{lem:Gorenstein_projective_epic_1} equal to
$\kercomplx^i \stackrel{k^i}{\longrightarrow} \Image
\partial_{\kercomplx}^i$.  So $k^i$ is a relative epimorphism with
respect to the class of Gorenstein projective modules.

The one short exact sequence which is missing is
\[
  0 \rightarrow \Ker u
    \longrightarrow \Ker \partial_{\kercomplx}^1
    \stackrel{u}{\longrightarrow} \M
    \rightarrow 0.
\]
Here lemma \ref{lem:zero} gives that lemma
\ref{lem:Gorenstein_projective_epic_2} applies if I set the
homomorphism $\Ker \partial_{\kercomplx}^i
\stackrel{t}{\longrightarrow} \N$ from lemma
\ref{lem:Gorenstein_projective_epic_2} equal to $\Ker 
\partial_{\kercomplx}^1 \stackrel{u}{\longrightarrow} \M$.  So $u$ is
a relative epimorphism with respect to the class of Gorenstein
projective modules.
\end{proof}

\begin{Remark}
\label{rmk:special_precovers}
For the next theorem, recall that a Gorenstein projective precover of
an \Rlm\ $\M$ is a homomorphism $\Gorprojmod \longrightarrow \M$ which
is a relative epimorphism with respect to the class of Gorenstein
projective modules, where $\Gorprojmod$ is Gorenstein projective
(cf.\ remark \ref{rmk:Gorenstein_projective_epic}). 

In \cite[thm.\ 3.1 and lem.\ 4.1]{AvrMart} and
\cite[thm.\ (2.10)]{HHGorensteinHomDim} is proved the following
``approximation'' result: Each module with finite Gorenstein
projective dimension has a Gorenstein projective precover whose kernel
has finite projective dimension.  

On the other hand, the existence of such a precover forces a module to
have finite Gorenstein projective dimension, so such precovers cannot
exist more generally.

However, one salient feature of modules of finite projective dimension
is that their projective resolutions are also Gorenstein projective
resolutions.  So an alternative view of \cite[thm.\ 3.1 and lem.\
4.1]{AvrMart} and \cite[thm.\ (2.10)]{HHGorensteinHomDim} is that they
prove that each module with finite Gorenstein projective dimension has
a Gorenstein projective precover whose kernel is nice, namely, each
projective resolution of the kernel is also a Gorenstein projective
resolution.

This result generalizes to all modules as follows.  Note again that
some authors refer to my Gorenstein projective resolutions as
``proper'', cf.\ remark \ref{rmk:Gorenstein_projective_resolutions}. 
\end{Remark}

\begin{Theorem}
\label{thm:special_precovers}
Let $\M$ be an \Rlm.  Then there is a Gorenstein projective
precover $\Gorprojmod \longrightarrow \M$ whose kernel $\kermod$
satisfies that each projective resolution of $\kermod$ is also a
Gorenstein projective resolution of $\kermod$.
\end{Theorem}

\begin{proof}
Since all projective resolutions of $\kermod$ are isomorphic in
$\KProjR$, it is not hard to see that it is sufficient to produce a
Gorenstein projective precover whose kernel $\kermod$ has {\em one}
projective resolution which is also a Gorenstein projective
resolution.

However, the exact sequence from remark \ref{rmk:resolution} and lemma
\ref{lem:resolution} splits into exact sequences
\[
  0 \rightarrow \Ker u
    \longrightarrow \Ker \partial_{\kercomplx}^1
    \stackrel{u}{\longrightarrow} \M
    \rightarrow 0
\]
and
\begin{equation}
\label{equ:resolution2}
  \cdots 
  \longrightarrow \kercomplx^{-2}
  \longrightarrow \kercomplx^{-1}
  \longrightarrow \kercomplx^0
  \longrightarrow \Ker u
  \rightarrow 0.
\end{equation}
Lemma \ref{lem:resolution} implies that $\Ker \partial_{\kercomplx}^1
\stackrel{u}{\longrightarrow} \M$ is a Gorenstein projective
precover, and that the sequence \eqref{equ:resolution2} is an
augmented Gorenstein projective resolution of the kernel $\Ker u$.
But \eqref{equ:resolution2} is clearly also an augmented projective
resolution of $\Ker u$ because the modules $\kercomplx^i$ are
projective, and so the theorem follows with $\Gorprojmod
\longrightarrow \M$ equal to $\Ker \partial_{\kercomplx}^1
\stackrel{u}{\longrightarrow} \M$ and $\kermod$ equal to $\Ker
u$.
\end{proof}

\section{Comparing absolute, relative, and Tate cohomology}
\label{sec:sequence}

This section proves that the new definition of Tate $\Ext$ groups of
this paper fits into the long exact sequence
\eqref{equ:long_exact_sequence} from the introduction (theorem
\ref{thm:long_exact_sequence}).

\begin{Remark}
In this section appear the relative $\Ext$ functors $\GExt$.  They are
defined by
\[
  \GExt^i(\M,-) = \H^i \Hom_{\R}(\Gorprojcomplx,-)
\]
where $\Gorprojcomplx$ is a Gorenstein projective resolution of $\M$
(cf.\ remark \ref{rmk:Gorenstein_projective_resolutions}).  Such a
resolution exists by lemma \ref{lem:resolution}.  See
\cite{AvrMart} or \cite{HHGorensteinDerived} for an exposition of the
theory of $\GExt^i$.
\end{Remark}

\begin{Construction}
\label{con:truncated_sequence}
Consider the short exact sequence from construction
\ref{con:basic_exact_sequence},
\[
  0 \rightarrow \kercomplx 
    \longrightarrow \Enochscomplxtil
    \longrightarrow \projcomplx \rightarrow 0.
\]
Truncating the complexes $\kercomplx$ and $\Enochscomplxtil$ gives a new
short exact sequence of complexes,
\[
  \begin{diagram}[labelstyle=\scriptstyle,height=4ex,width=7ex]
    & & \vdots & & \vdots & & \vdots & & \\
    & & \uTo & & \uTo & & \uTo & & \\
    0 & \rTo & 0 & \rTo & 0 & \rTo & 0 & \rTo & 0 \\
    & & \uTo & & \uTo & & \uTo & & \\
    0 & \rTo & \Ker \partial_{\kercomplx}^1 & \rTo & \Ker \partial_{\Enochscomplxtil}^1 & \rTo & 0 & \rTo & 0 \\
    & & \uTo & & \uTo & & \uTo & & \\
    0 & \rTo & \kercomplx^0 & \rTo & \Enochscomplxtil^0 & \rTo & \projcomplx^0 & \rTo & 0 \\
    & & \uTo & & \uTo & & \uTo & & \\
    0 & \rTo & \kercomplx^{-1} & \rTo & \Enochscomplxtil^{-1} & \rTo & \projcomplx^{-1} & \rTo & 0 \lefteqn{,} \\
    & & \uTo & & \uTo & & \uTo & & \\
    & & \vdots & & \vdots & & \vdots & & \\
  \end{diagram}
\]
which I will denote
\begin{equation}
\label{equ:truncated_sequence}
  0 \rightarrow \kercomplxprime \longrightarrow \Enochscomplxtilprime
  \longrightarrow \projcomplx \rightarrow 0.
\end{equation}
\end{Construction}

\begin{Theorem}
\label{thm:long_exact_sequence}
Let $\M$ and $\N$ be \Rlm s.  Then there is a long exact
sequence 
\begin{eqnarray*}
  0 & \longrightarrow & \GExt^1(\M,\N) 
                        \longrightarrow \Ext^1(\M,\N)
                        \longrightarrow \TExt^1(\M,\N) \\
    & \longrightarrow & \cdots \\
    & \longrightarrow & \GExt^i(\M,\N) 
                        \longrightarrow \Ext^i(\M,\N)
                        \longrightarrow \TExt^i(\M,\N) 
                        \longrightarrow \cdots, \\
\end{eqnarray*}
natural in $\M$ and $\N$.
\end{Theorem}

\begin{proof}
Consider the short exact sequence \eqref{equ:truncated_sequence} from
construction \ref{con:truncated_sequence}.  Let me set
\[
  \projcomplx = \res \M, 
\]
where $\res \M$ is the projective resolution depending functorially on
$\M$.  Recall that $\Enochscomplxtilprime$ is a truncation of
\[
  \Enochscomplxtil = \inc \adj \projcomplx = \adj \res \M.
\]
Since $\res \M$ consists of projective modules, the short exact
sequence \eqref{equ:truncated_sequence} is semi-split and therefore
stays exact under the functor $\Hom_{\R}(-,\N)$.  So there is a short
exact sequence of complexes
\begin{equation}
\label{equ:Hom_of_truncated_sequence}
  0 \rightarrow \Hom_{\R}(\res \M,\N)
    \longrightarrow \Hom_{\R}(\Enochscomplxtilprime,\N)
    \longrightarrow \Hom_{\R}(\kercomplxprime,\N)
    \rightarrow 0.
\end{equation}

Since $\res \M$ is a projective resolution of $\M$, I have
\[
  \H^i\!\Hom_{\R}(\res \M,\N) = \Ext^i(\M,\N)
\]
for each $i$.  Moreover, the form of the complex $\Enochscomplxtilprime$
makes it clear that $\H^0\!\Hom_{\R}(\Enochscomplxtilprime,\N) = 0$
and that 
\begin{eqnarray*}
  \H^i\!\Hom_{\R}(\Enochscomplxtilprime,\N)
  & = & \H^i\!\Hom_{\R}(\Enochscomplxtil,\N) \\
  & = & \H^i\!\Hom_{\R}(\adj \res \M,\N) \\
  & = & \TExt^i(\M,\N)
\end{eqnarray*}
for $i \geq 1$.

Finally, lemma \ref{lem:resolution} says that $\kercomplxprime$ is a
Gorenstein projective resolution of $\M$, shifted by one.  Hence
\[
  \H^i\!\Hom_{\R}(\kercomplxprime,\N)
  = \GExt^{i+1}(\M,\N)
\]
for $i \geq -1$.  So looking at the cohomology long exact sequence of
\eqref{equ:Hom_of_truncated_sequence}, starting with
$\H^0\!\Hom_{\R}(\Enochscomplxtilprime,\N) = 0$, gives
\[
  0 \rightarrow \GExt^1(\M,\N) 
    \longrightarrow \Ext^1(\M,\N)
    \longrightarrow \TExt^1(\M,\N)
    \longrightarrow \cdots,
\]
that is, the sequence from the theorem.
\end{proof}

\bigskip

\noindent
{\bf Acknowledgement.}  
The diagrams were typeset with Paul Taylor's
{\tt diagrams.tex}.

\end{document}